\documentclass[a4paper,11pt]{article}
\usepackage{amsmath,amssymb,amsthm,psfrag,color,graphicx,caption2}

\renewcommand{\captionlabeldelim}{~}

\newtheorem*{theoremA}{Theorem A}
\newtheorem*{theoremB}{Theorem B}
\newtheorem*{theoremC}{Theorem C}
\newtheorem*{theoremD}{Gr\"{o}tzsch Theorem}
\newtheorem*{theorem1}{Theorem 1}
\newtheorem*{theorem2}{Theorem 2}
\newtheorem*{conj1}{Conjecture 1}
\newtheorem*{conj2}{Conjecture 2}
\newtheorem*{conj3}{Conjecture 3}

\newtheorem{lemma}{Lemma}
\newtheorem{cor}{Corollary}
\newtheorem{prop}{Proposition}

\theoremstyle{definition}
\newtheorem{definition}{Definition}
\newtheorem*{remark1}{Remark 1}
\newtheorem*{remark2}{Remark 2}

\newtheorem*{acknowledgements}{Acknowledgements}

\numberwithin{equation}{section}

\begin{document}

\title{No invariant line fields on Cantor Julia sets}

\author{Yongcheng Yin\and Yu Zhai}
\date{}
\maketitle

\begin{abstract}
In this paper, we prove that a rational map with a Cantor Julia
set carries no invariant line field on its Julia set. It follows
that a structurally stable rational map with a Cantor Julia set is
hyperbolic.
\end{abstract}

\vspace*{+3mm}

\section{Introduction and statements}\label{intro}

\noindent Let $f:\hat{\mathbb{C}}\rightarrow \hat{\mathbb{C}}$ be
a rational map of the Riemann sphere to itself of degree
$d\geqslant 2$. The map $f$ is hyperbolic if there are a smooth
conformal metric $\rho$ defined on a neighborhood of the Julia set
$J(f)$ and a constant $C>1$ such that $\|f'(z)\|_{\rho}>C$ for all
$z\in J(f)$. It is equivalent to every critical point of $f$ tends
to an attracting periodic cycle under forward iteration. See
\cite{Mc1} and \cite{Mi}.

Let $\mathrm{Rat}_{d}$ be the space of all the M\"{o}bius
equivalence classes of rational maps of degree $d$. The space
$\mathrm{Rat}_{d}$ has dimension $2d-2$. A central problem in
holomorphic dynamics is the following.

\begin{conj1}[Density of hyperbolicity] The set of
hyperbolic rational maps is open and dense in the space
$\mathrm{Rat}_{d}$.
\end{conj1}

Openness of the set of hyperbolic rational maps is known, but
density is only known in the family of real polynomials, see
\cite{GSw}, \cite{KSS1}, \cite{KSS2}, \cite{Ly1}, \cite{Mc1} and
\cite{Shi}.

A rational map $f$ admits an {\it{invariant line field}} on the
Julia set $J(f)$ if there is a measurable Beltrami differential
$\mu (z)\frac{d\bar{z}}{dz}$ on $\hat{\mathbb{C}}$ such that
$f^*\mu =\mu$ a.e., $|\mu|=1$ on a positive measure subset $E$ of
$J(f)$ and $\mu =0$ on $\hat{\mathbb{C}}\setminus E$.

A rational map $f$ is called a {\it{Latt$\grave{e}$s example}} if
it is doubly covered by an integral torus endmorphism. The Julia
set of such a rational map is $\hat{\mathbb{C}}$ and
$\frac{d\bar{z}}{dz}$ is an invariant line field on
$\hat{\mathbb{C}}$.

\begin{conj2} A rational map $f$ carries no invariant line
fields on its Julia set, except when $f$ is a Latt$\grave{e}$s
example.
\end{conj2}

This conjecture is stronger than the density of hyperbolic
dynamics.

\begin{theoremA}[\cite{McS}] The no invariant line field conjecture
implies the density of hyperbolic dynamics in the space of all
rational maps.
\end{theoremA}

The absence of invariant line fields on the Julia set is known in
the following cases:
\begin{enumerate}
\item[(1)] Non-infinitely renormalizable quadratic polynomials with no
irrational indifferent periodic points, \cite{Ly2}, \cite{Shi} and
\cite{Yo}.

\item[(2)] Robust infinitely renormalizable quadratic polynomials
and real quadratic polynomials, \cite{Mc1}.

\item[(3)] Quadratic polynomials with a Siegel cycle of bounded
type rotation number, \cite{LY}, \cite{Mc3} and \cite{Pe}.

\item[(4)] Real polynomials with only one non-escaping critical
point which is real and has odd local degree, \cite{LS}.

\item[(5)] Real rational maps(non Latt$\grave{e}$s example)
whose critical points are all on the extended real axis and have
even local degrees, \cite{Sh}.

\item[(6)] Summable rational maps with completely invariant Fatou
domains, \cite{Mak1}.

\item[(7)] Summable rational maps with small exponents,
\cite{GSm}.

\item[(8)] Weakly hyperbolic rational maps, \cite{Hai}.
\end{enumerate}

In this paper we will prove

\begin{theorem1}
Let $f$ be a rational map with a Cantor Julia set. Then $f$
carries no invariant line fields on its Julia set.
\end{theorem1}

\begin{remark1}
(1) It is reasonable to conjecture that a Cantor Julia set always
has measure zero. Theorem 1 can be regarded as a step towards this
conjecture.

(2) As a special case of Theorem 1, we disprove a question in the
Chapter 12 of \cite{BH}.
\end{remark1}

Let $X$ be a complex manifold. A {\it{holomorphic family}} of
rational maps $\{f_{\lambda}(z)\}_{\lambda \in X}$ is a
holomorphic map $X\times \hat{\mathbb{C}}\rightarrow
\hat{\mathbb{C}}$, given by $(\lambda, \,z)\rightarrowtail
f_{\lambda}(z)$. Let $X^s\subset X$ be the set of structurally
stable parameters. That is, $a\in X^s$ if and only if there is a
neighborhood $U$ of $a$ such that $f_a$ and $f_b$ are
topologically conjugate for all $b\in U$. The space $X^{qc}\subset
X$ of quasiconformally stable parameters is defined similarly,
with conjugacy replaced by quasiconformality.

Using the Harmonic $\lambda$-Lemma of Bers and Royden, McMullen
and Sullivan proved the following result.

\begin{theoremB}[\cite{MSS} and \cite{McS}]
In any holomorphic family of rational maps, $X^s$ is open and
dense in $X$. Moreover, the structurally stable and
quasiconformally stable parameters coincide, i.e., $X^s=X^{qc}$.
\end{theoremB}

\begin{conj3} A structurally stable rational map is
hyperbolic.
\end{conj3}

Combine Theorem 1 and the theory of Teichm\"{u}ller space of
rational maps in \cite{McS}, we can prove the following result.

\begin{theorem2}
Let $f$ be a rational map with a Cantor Julia set. If $f$ is
structurally stable, then it is hyperbolic.
\end{theorem2}

\begin{remark2} The same result in Theorem 2 was also proved by
Makienko under some additional assumptions, see \cite{Mak2}.
\end{remark2}

A nested sequence of some critical pieces constructed by
Kozlovski, Shen, and van Strien in \cite{KSS2}, which we shall
call ``KSS nest", will play a crucial rule in the proof of Theorem
1. Principal nest and modified principal nest are used to study
the dynamics of unicritical polynomials, see \cite{AKLS},
\cite{BSS}, \cite{KL} and \cite{Ly1}. In \cite{Ly1}, Lyubich
proved the linear growth of its ``principal moduli" for quadratic
polynomials. This yields the density of hyperbolic maps in the
real quadratic family. The same result is also obtained by Graceyk
and \'{S}wi\c{a}tek in \cite{GSw}. See also \cite{Mc1} and
\cite{Shi}. Recently, the local connectivity of Julia sets and
combinatorial rigidity for unicritical polynomials are proved in
\cite{KL} and \cite{AKLS} by means of principal nest and modified
principal nest.

This paper is organized as follows. In section 2, we present some
distortion lemmas which are used in section 4. In section 3, we
introduce the Branner-Hubbard puzzle about rational maps with
Cantor Julia set and the KSS nest constructed in \cite{KSS2}. By
means of the KSS nest and the distortion lemmas, we prove that the
shapes of some critical puzzle pieces are bounded in section 4. In
section 5, we give the proofs of Theorem 1 and Theorem 2.

\section{Distortion lemmas }

\noindent Any doubly connected domain $A$ on the complex plane is
conformally equivalent to one of the following three types of
typical domains:
\begin{enumerate}
\item[(1)] $\mathbb{C}\backslash \{0\}$,

\item[(2)] $\Delta\backslash \{0\}$, where
$\Delta=\{z\,|\,|z|<1\}$,

\item[(3)] $A_{R}=\{z\,|\,1<|z|<R\}$.
\end{enumerate}

In the case $A$ is conformally equivalent to $A_R$, the modulus of
$A$ is defined as $\mathrm{mod}(A)=\frac{1}{2\pi}\ln R$. In the
other two cases, $\mathrm{mod}(A)= \infty$.

For $0<r<1$, let $B_{r}=\Delta \setminus [0,r]$. The modulus
$\mathrm{mod}(B_r)$ is decreasing in $(0,1)$ with $\lim_{r\to
1^-}\mathrm{mod}(B_r)=0$ and $\lim_{r\to
0^+}\mathrm{mod}(B_r)=+\infty$.

\begin{theoremD}[\cite{Ah}] Let $A$ be a doubly connected domain in
$\Delta$ which separates the unit circle from the points
$\{0,r\}$. Then $\mathrm{mod}(A)\leqslant \mathrm{mod}(B_r)$.
\end{theoremD}

Denote $d(\cdot,\cdot)$ and $\mathrm{diam}$ the distance and
diameter with respect to the Euclidean metric respectively.

\begin{lemma}
Let $\widetilde{U}\subset\subset U\neq \mathbb{C}$ be two simply
connected domains and $\mathrm{mod}(U\setminus
\overline{\widetilde{U}})\geqslant m>0$. Then there exists a
constant $c=c(m)>0$ such that
$$d(\omega,\partial U)\geqslant c\cdot \mathrm{diam}(\widetilde{U})$$
for any $\omega\in \overline{\widetilde{U}}$.
\end{lemma}

\begin{proof}
For any $\omega\in \overline{\widetilde{U}}$, let $h_{\omega}(z)$
be a conformal map from $\Delta$ onto $U$ with
$h_{\omega}(0)=\omega$. Then
$$\mathrm{mod}(\Delta\setminus
\overline{h_{\omega}^{-1}(\widetilde{U})})=\mathrm{mod}(U\setminus
\overline{\widetilde{U}})\geqslant m>0.$$ From Gr\"{o}tzsch
Theorem, there exists a constant $r_0=r_0(m)<1$ such that
$h_{\omega}^{-1}(\widetilde{U})\subset \{z\,|\,|z|<r_0\}$. By
Koebe Distortion Theorem, we have
$$d(\omega,\partial U)\geqslant
\frac{1}{4}|h_{\omega}'(0)|$$ and
$$\mathrm{diam}(\widetilde{U})\leqslant
2|h_{\omega}'(0)|\frac{r_0}{(1-r_0)^{2}}.$$

So we can take $\displaystyle c=c(m)=\frac{(1-r_0)^{2}}{8r_0}>0$,
which satisfies the inequality
$$d(\omega,\partial U)\geqslant c\cdot \mathrm{diam}(V).$$
\end{proof}

Let $U$ be a simply connected domain and $\omega\in U$. The {\it
Shape} of $U$ about $\omega$, denoted by
$\mathrm{Shape}(U,\omega)$, is defined as
$$\mathrm{Shape}(U,\omega)=\frac{\max_{z\in
\partial U} d(\omega,z)}{\min_{z\in
\partial U} d(\omega,z)}=\frac{\max_{z\in \partial U} d(\omega,z)}{d(\omega,\partial U)}.$$

\begin{lemma} Let $g:(\Delta,U,\widetilde{U})\rightarrow (\Delta,V,\widetilde{V})$
be a holomorphic proper map of degree $d$ with
$0\in\widetilde{U}\subset U \subset \Delta$, $0\in
\widetilde{V}\subset V \subset \Delta$. Suppose that
\begin{enumerate}
\item[(1)] $\deg(g\mid_{\widetilde{U}})=\deg(g\mid_{U})=\deg(g\mid_{\Delta})=d\geqslant 2$,

\item[(2)] $\mathrm{mod}(V\setminus \overline{\widetilde{V}})\geqslant m>0$.
\end{enumerate}
Then there exists a constant $K=K(m,d)>0$ such that
$$\mathrm{Shape}(U,0)\leqslant K\cdot
\mathrm{Shape}(V,0)^{\frac{1}{d}}.$$
\end{lemma}

\begin{proof}
Let
\begin{align*}
R=\max_{z\in\partial U}{|z|}&,\textrm{ }r=\min_{z\in \partial U}{|z|};\\
R'=\max_{\omega\in \partial
V}{|\omega|}&,\textrm{}r'=\min_{\omega\in
\partial V}{|\omega|}.
\end{align*}
There are points $z_{R}\in
\partial U$ and $z_{r}\in \partial
U$ such that $R=|z_{R}|$ and $r=|z_{r}|$.

The holomorphic proper map $g$ can be written as
$$g(z)=e^{i\theta}\prod_{j=1}^{d}\frac{z-a_{j}}{1-\overline{a}_{j}z},$$
where $\theta\in [0,2\pi)$ and $a_{j}\in
\Delta,\,j=1,\cdot\cdot\cdot,d.$ By the first assumption, we have
$a_{j}\in \widetilde{U}$ and $|a_{j}|< |z_{R}|<1$ for all $j$.

Consider the annulus $\{z\,|\,|z|<R\}\setminus
\overline{\widetilde{U}}$, we have
$$\mathrm{mod}(\{z\,|\,|z|<R\}\setminus
\overline{\widetilde{U}})\geqslant \mathrm{mod}(U\setminus
\overline{\widetilde{U}})\geqslant \frac{m}{d}.$$ By Gr\"{o}tzsch
Theorem, there exists a constant $c_{0}=c_{0}(\frac{m}{d})>1$ such
that $R=|z_{R}|\geqslant c_{0}|a_{j}|$ for all $j$. We have
\begin{align*}
R'&\geqslant|g(z_{R})|=\prod_{j=1}^{d}\frac{|z_{R}-a_{j}|}{|1-\overline{a}_{j}z_{R}|}\geqslant
\prod_{j=1}^{d}\frac{|z_{R}|-|a_{j}|}{1+|a_{j}z_{R}|}\\
&\geqslant\prod_{j=1}^{d}\frac{|z_{R}|-|a_{j}|}{2}\geqslant\prod_{j=1}^{d}(\frac{c_{0}-1}{2c_{0}})|z_{R}|=c_{1}^{d}R^{d}
\end{align*}
with $\displaystyle c_{1}=c_{1}(m,d)=\frac{c_{0}-1}{2c_{0}}$. On
the other hand,
\begin{align*}
r'&\leqslant|g(z_{r})|=\prod_{j=1}^{d}\frac{|z_{r}-a_{j}|}{|1-\overline{a}_{j}z_{r}|}\leqslant
\prod_{j=1}^{d}\frac{r+|a_j|}{1-|a_j|}\leqslant
\prod_{j=1}^{d}\frac{r+\frac{r}{c}}{1-r_0}=c_{2}^{d}r^{d}
\end{align*}
with $\displaystyle c_{2}=c_{2}(m,d)=\frac{c+1}{c(1-r_{0})}$ in
which $c$ and $r_0$ come from Lemma 1 and its proof. It follows
that
$$\mathrm{Shape}(U,0)\leqslant K\cdot \mathrm{Shape}(V,0)^{\frac{1}{d}}$$
with $\displaystyle K=K(m,d)=\frac{c_{2}}{c_{1}}$.
\end{proof}

\section{Branner-Hubbard puzzle and KSS nest}

\noindent From now on, we always assume that the Julia set $J(f)$
of a rational map $f$ is a cantor set. The Fatou set $F(f)$ has
only one component. It is either an attracting basin or a
parabolic basin.

We first construct the Branner-Hubbard puzzle.

{\bf{The attracting case}}. We assume that $\infty$ is the fixed
attracting point. Take a simply connected neighborhood
$U_{0}\subset F(f)$ of $\infty$ such that $U_{0}\subset\subset
f^{-1}(U_{0})$. Let $U_{n}$ be the component of $f^{-n}(U_{0})$
containing $\infty$. Then $U_{n}\subset\subset U_{n+1}$ and
$\displaystyle F(f)=\cup_{n=0}^{\infty}U_{n}$. For a large enough
integer $N_{0}$, $f^{-n}(U_{N_{0}})$ has only one component for
any $n\geqslant 0$. The set $f^{-n}(\hat{\mathbb{C}}\setminus
\overline{U}_{N_{0}})$ is the disjoint union of a finite number of
topological disks. For each $n\geqslant 0$, let $\mathbf{P_{n}}$
be the collection of all components of
$f^{-n}(\hat{\mathbb{C}}\setminus \overline{U}_{N_{0}})$ which are
called puzzle pieces of depth $n$.

For any point $x\in J(f)$ and any $n\geqslant 0$, there is only
one $P_{n}(x)\in \mathbf{P_{n}}$ containing $x$. Thus each point
$x\in J(f)$ determines a nested sequence $P_{0}(x)\supset
P_{1}(x)\supset \cdot\cdot\cdot$ and $\cap_{n\geqslant
0}P_n(x)=\{x\}$.

{\bf{The parabolic case}}. We suppose that $0$ is the parabolic
fixed point and $\infty$ is in the Fatou set. According to the
Leau-Fatou Flower Theorem, there is a flower petal $U_{0}\subset
F(f)$ with $0\in \partial U_{0}$ such that
$\overline{U_{0}}\subset f^{-1}(U_{0})\cup \{0\}$. We can
construct the puzzle as in the attracting case. Each point $x\in
J(f)\setminus \cup_{n\geqslant 0}f^{-n}(0)$ determines a nested
sequence $P_{0}(x)\supset P_{1}(x)\supset \cdot\cdot\cdot$ and
$\cap_{n\geqslant 0}P_n(x)=\{x\}$.

Take $N_0$ large enough such that $U_{N_{0}}$ contains all
critical points in the Fatou set and each puzzle piece contains at
most one critical point.

Let
$$\mathrm{Crit}=\{c\in J(f)\,|\,c\textrm{ is the critical point of }f\}$$
in the attracting case and
$$\mathrm{Crit}=\{c\in J(f)\setminus \cup_{n\geqslant 0}f^{-n}(0)\,|\,c\textrm{ is the critical point of }f\}$$
in the parabolic case.

For each $x\in J(f)$ (in parabolic case, $\displaystyle x\in
J(f)\setminus \cup_{n\geqslant 0}f^{-n}(0)$ respectively), the
tableaux $T(x)$ is defined in \cite{BH}. It is a
 two dimension
array $\{P_{n,l}(x)\}_{n\geqslant0,l\geqslant 0}$ with
$P_{n,l}(x)=f^{l}(P_{n+l}(x))=P_n(f^l(x))$. The position $(n,l)$
is called critical if $P_{n,l}(x)$ contains a critical point of
$f$. If $P_{n,l}(x)$ contains a critical point $c$, the position
$(n,l)$ is called a $c$-position. The tableau $T(c)$ of a critical
point $c\in \mathrm{Crit}$ is called periodic if there is a
positive integer $k$ such that $P_{n}(c)=f^{k}(P_{n+k}(c))$ for
all $n\geqslant 0$. Since the Julia set is a Cantor set, $T(c)$ is
not periodic for all $c\in \mathrm{Crit}$.

All the tableaus satisfy the following three rules
\begin{enumerate}
\item[(T1)] If $P_{n,l}(x)=P_n(c)$ for some critical point $c$,
then $P_{i,l}(x)=P_i(c)$ for all $0\leqslant i \leqslant n$.

\item[(T2)] If $P_{n,l}(x)=P_n(c)$ for some critical point $c$,
then $P_{i,l+j}(x)=P_{i,j}(c)$ for $i+j\leqslant n$.

\item[(T3)] Let $T(c)$ be a tableau for some critical point $c$
and $T(x)$ be any tableau. Assume
           \begin{enumerate}
           \item $P_{n+1-l,l}(c)=P_{n+1-l}(c_1)$ for some critical
                 point $c_1$ and $n>l\geqslant
                 0$, and $P_{n-i,i}(c)$ contains no critical points for
                 $0<i<l$.
           \item $P_{n,m}(x)=P_n(c)$ and $P_{n+1,m}(x)\neq
                  P_{n+1}(c)$ for some $m>0$.
           \end{enumerate}
Then $P_{n+1-l,m+l}(x)\neq P_{n+1-l}(c_1)$.
\end{enumerate}

\begin{definition}
(1)\; The tableau $T(x)$ for $x$ is {\it{non-critical}} if there
exists an integer $n_0\geqslant 0$ such that $(n_0,j)$ is not
critical for all $j>0$.

(2)\; We write $x\to y$ if for any $n\geqslant 0$, there exists
$j>0$ such that $y\in P_{n,j}(x)$, i.e., $f^j(P_{n+j}(x))=
P_n(y)$. It is clear that $x\to y$ if and only if $y\in
\cup_{n>0}f^{-n}(x)$ or $y\in \omega (x)$, the limit set of the
forward orbit of $x$. If $x\to y$ and $y\to z$, then $x\to z$. For
each critical point $c\in\mathrm{Crit}$, let
$$F(c)=\{c^\prime\in\mathrm{Crit}\,| \,\, c\to c^\prime\}$$
and
$$[c]=\{c^\prime\in\mathrm{Crit}\,| \,\, c\to c^\prime \textrm{ and } c^\prime \to
c \}.$$

(3)\;We say $P_{n+k}(c^\prime)$ is a {\it{child}} of $P_n(c)$ if
$c^\prime\in [c]$, $f^k(P_{n+k}(c^\prime))= P_n(c)$, and
$f^{k-1}:\, P_{n+k-1}(f(c^\prime))\to P_n(c)$ is conformal.

(4)\; Suppose $c\to c$ , i.e., $[c]\not=\emptyset$. We say $T(c)$
is {\it{persistently recurrent}} if $P_n(c_1)$ has only finitely
many children for all $n\geqslant 0$ and all $c_1\in [c]$.
Otherwise, $T(c)$ is said to be {\it{reluctantly recurrent}}.
\end{definition}

Take $N_{0}$ large enough such that for any $c\in \mathrm{Crit}$,
there is no $c'$-position in the first row of $T(c)$ if $c\not\to
c'$.

Let
\begin{eqnarray*}
&&\textrm{Crit}_\textrm{n}=\{c\in\mathrm{Crit}\,| \,\, T(c) \textrm { is non-critical}\},\\
&&\textrm{Crit}_\textrm{p}=\{c\in\mathrm{Crit}\,| \,\, T(c) \textrm { is persistently recurrent}\},\\
&&\textrm{Crit}_\textrm{r}=\{c\in\mathrm{Crit}\,| \,\, T(c) \textrm { is reluctantly recurrent}\},\\
&&\textrm{Crit}_\textrm{en}=\{c^\prime\in\mathrm{Crit}\,|\, \,
c^\prime \not\to c^\prime
\textrm{ and } c^\prime \to c \textrm{ for some } c\in\textrm{Crit}_\textrm{n} \},\\
&&\textrm{Crit}_\textrm{ep}=\{c^\prime\in\mathrm{Crit}\,| \,\,
c^\prime \not\to c^\prime
\textrm{ and } c^\prime \to c \textrm{ for some } c\in \textrm{Crit}_\textrm{p}\},\\
&&\textrm{Crit}_\textrm{er}=\{c^\prime\in\mathrm{Crit}\,| \,\,
c^\prime \not\to c^\prime \textrm{ and } c^\prime \to c \textrm{
for some } c\in \textrm{Crit}_\textrm{r}\}.
\end{eqnarray*}
Then
$$\mathrm{Crit}=\textrm{Crit}_\textrm{n}\cup \textrm{Crit}_\textrm{p}
\cup \textrm{Crit}_\textrm{r}\cup \textrm{Crit}_\textrm{en} \cup
\textrm{Crit}_\textrm{ep}\cup \textrm{Crit}_\textrm{er}.$$ This is
not a classification because these sets might intersect.

The following lemma can be found in \cite{QY}.

\begin{lemma}
If $T(c)$ is persistently recurrent, then  $F(c)=[c]$.
\end{lemma}

Now, we briefly introduce a critical nest which is constructed by
Kozlovski, Shen, and van Strien in \cite{KSS2}. Such nest will be
called {\it KSS} nest.

Let $A$ be an open set and $x\in A$. The connected component of
$A$ containing $x$ will be denoted by Comp$_{x}(A)$. Given a
puzzle piece $I$, let
$$D(I)=\{z\in \mathbb{C}\,|\,\exists k\geqslant
1\textrm{ }s.t.\textrm{ }f^{k}(z)\in I\}=\bigcup_{k\geqslant
1}f^{-k}(I).$$ For any $z\in D(I)$, let $\mathcal{L}_{z}(I)$ be
the connected component of $D(I)$ containing $z$. We further
define $\hat{\mathcal{L}}_{z}(I)=I$ if $z\in I$ and
$\hat{\mathcal{L}}_{z}(I)=\mathcal{L}_{z}(I)$ if $z\in
D(I)\setminus I$.

For any $z\in D(I)$, let $k\geqslant 1$ be the integer such that
$f^k(\mathcal{L}_{z}(I))=I$ and let $n_0$ be the depth of $I$. By
tableau rules (T1) and (T2), there is at most one $c$-position on
the diagonal
$$\{(n,m)\,|\,\, n+m = n_0+k, \quad n_0< n \leqslant n_0+k\}$$
in the tableau $T(z)$ for any $c\in \textrm{Crit}$. Hence
$$\deg(f^k:\mathcal{L}_{z}(I)\rightarrow I)\leqslant D$$
for some constant $D<\infty$ depending only on $\textrm{Crit}$.

Suppose $T(c_0)$ is persistently recurrent, then $F(c_0)=[c_0]$.
Let
$$b=\#[c_0],\,\,\,\,d_0=\deg_{c_0}f,\,\,\,\,d_{max}=\max\{\deg_cf
|\,c\in [c_0]\},$$ and
$$\textrm{orb}([c_0])=\bigcup_{n\geqslant 0}f^n([c_0]).$$

For each puzzle piece $I\ni c_0$, there are pullbacks
$\mathcal{A}(I)\subset \mathcal{B}(I)$ of $I$ containing $c_0$
with the following properties
\begin{enumerate}

\vspace{-0.2cm} \item[(P1)]\; $f^t (\mathcal{B}(I))= I$ and
$\deg(f^t|_{\mathcal{B}(I)}) \leqslant {d_{max}}^{b^2}$,

\vspace{-0.2cm}\item[(P2)]\;
$\mathcal{A}(I)=\textrm{Comp}_{c_0}f^{-t}(\mathcal{L}_{f^t(c_0)}(I))$,
$f^s(\mathcal{A}(I))= I$ , $s-t\geqslant 1$ is the smallest
integer such that $f^{s-t}(f^t(c_0))\in I$ and
$\deg(f^s|_{\mathcal{A}(I)}) \leqslant {d_{max}}^{b^2+b}$,

\vspace{-0.2cm}\item[(P3)]\; $(\mathcal{B}(I)-\mathcal{A}(I)) \cap
\textrm{orb}([c_0]) =\emptyset$.
\end{enumerate}
For details, see \cite{KSS2} or \cite{QY}.

\begin{definition}
Given  a puzzle piece $P$ containing $c_0$, a {\it{successor}} of
$P$ is a piece of the form $\hat{\mathcal{L}}_{c_0}(Q)$, where $Q$
is a child of $\hat{\mathcal{L}}_c(P)$ for some $c\in [c_0]$.
\end{definition}

It is clear that $\mathcal{L}_{c_{0}}(P)$ is a successor of $P$.
Since $T(c_{0})$ is persistently recurrent, $P$ has at least two
successors and that $P$ has only finitely many successors. Let
$\Gamma(P)$ be the last successor of $P$. Then there exists an
integer $q\geqslant 1$, the largest among all of the successors of
$P$, such that $f^{q}(\Gamma(P))=P$. From the definition of
successor, we have

\begin{enumerate}

\vspace{-0.2cm} \item[(P4)]\;$\deg(f^{q}|_{\Gamma(P)})\leqslant
d_{max}^{2b-1}$.
\end{enumerate}

\begin{figure}[h]
\psfrag{In}{$I_n$}\psfrag{Tc0}{$T(c_0):$}

\psfrag{BIn}{$\mathcal{B}(I_n)$} \psfrag{fsn}{$f^{s_n}$}

\psfrag{LnAIn}{$L_n=\mathcal{A}(I_n)$} \psfrag{Knp}{$K_n'$}

\psfrag{ftn}{$f^{t_n}$} \psfrag{KnBLn}{$K_n=\mathcal{B}(L_n)$}

\psfrag{fbn}{$f^{b_n}$}
\psfrag{Knt}{$\widetilde{K}_n=\mathcal{A}(L_n)$}

\psfrag{gammaKn}{$\Gamma(K_n)$} \psfrag{fqn1}{$f^{q_{n,1}}$}

\psfrag{fqn}{$f^{q_n}$}
\psfrag{In1gammaTKn}{$I_{n+1}=\Gamma^{3b}(K_n)$} \centering
\includegraphics[width=10cm]{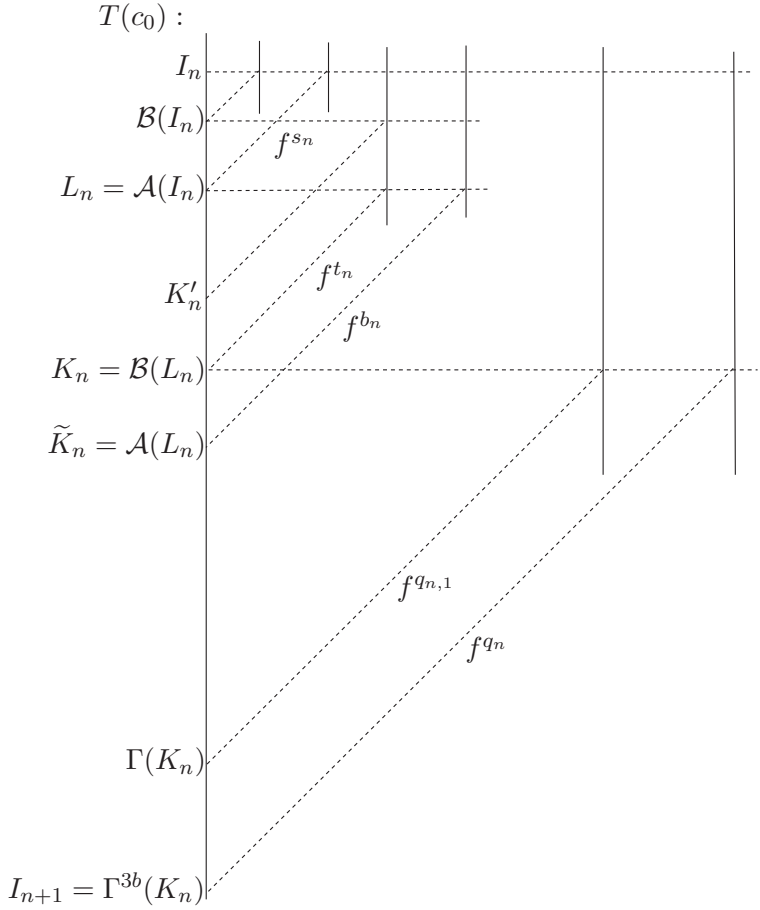} \caption{KSS
nest} \label{fig2}
\end{figure}

Now we can define the {\it KSS nest} in the following way: $I_{0}$
is a given piece containing $c_{0}$ and for $n\geqslant 0$,
\begin{eqnarray*}
&&L_{n}=\mathcal{A}(I_{n}),\\
&&M_{n,0}=K_{n}=\mathcal{B}(L_{n}),\\
&&M_{n,j+1}=\Gamma(M_{n,j}) \textrm{ for } 0\leqslant j\leqslant 3b-1,\\
&&I_{n+1}=M_{n,3b}=\Gamma^{3b}(K_{n})=\Gamma^{3b}(\mathcal{B}(\mathcal{A}(I_{n}))),
\end{eqnarray*}
with $b=\#[c_0]$. See Figure~\ref{fig2}.





Suppose $f^{s_{n}}(L_{n})=I_{n}$, $f^{t_{n}}(K_{n})=L_{n}$,
$f^{q_{n,j}}(M_{n,j})=M_{n,j-1}$ for $1\leqslant j\leqslant 3b$,
and $\displaystyle q_{n}=\sum_{j=1}^{3b}q_{n,j}$. Let
$K_{n}'=\mathrm{Comp}_{c_{0}}f^{-t_{n}}(\mathcal{B}(I_{n}))$ and
$\widetilde{K}_{n}=\mathcal{A}(L_{n})$. Then $(K'_n\setminus
K_n)\cap \mathrm{orb}([c_{0}])=\emptyset$ and
$$(K_{n}\setminus
\widetilde{K}_{n})\cap
\mathrm{orb}([c_{0}])=(\mathcal{B}(L_{n})\setminus
\mathcal{A}(L_{n}))\cap \mathrm{orb}([c_{0}])=\emptyset.$$ It
follows that
$$\deg(f^{t_{n}}|_{K_n'})=\deg(f^{t_n}|_{K_n}).$$
See Figure~\ref{fig2}.

\renewcommand{\captionlabeldelim}{~}

Let $p_{n}=q_{n-1}+s_{n}+t_{n}$. Then $f^{p_{n}}(K_{n})=K_{n-1}$
and
$$d_{0}^{3b+2}\leqslant \deg(f^{p_{n}}|_{K_{n}}) \leqslant
d_{1}=d_{max}^{8b^{2}-2b}.$$

For any puzzle piece $J$ containing $c_0$ and $z\in J\cap
\mathrm{orb}([c_{0}])$, let $r_z(J)=k(z)\geqslant 1$ be the
smallest integer such that $f^{k(z)}(z)\in J$ and
$$r(J)=\min\{k(z)\,|\,\,z\in J\cap
\mathrm{orb}([c_{0}])\}.$$ It is obvious that

\begin{enumerate}
\item[(1)] $r(J_1)\geqslant r(J_2)$ if $J_1\subset J_2$;
\item[(2)] $r(J)\geqslant k$ if $c_0\in J\subset J^\prime$, $f^k:\, J\to J^\prime$ and $c_0\not\in f^i (J)$ for
$0<i<k$;
\item[(3)] $\mathrm{depth}(\mathcal{A}(J))-\mathrm{depth}(\mathcal{B}(J))=r_{f^t(c_0)}(J)\geqslant
r(J)$.
\end{enumerate}

The following lemma plays a crucial rule in the proof of our
results and in \cite{QY}.

\begin{lemma}[\cite{KSS2}]
For any $n\geqslant 1$,
\begin{enumerate}
\item[(1)] $r(I_n)\leqslant s_n \leqslant (b+1)r(L_n)$;
\item[(2)] $r(L_n)\leqslant t_n \leqslant br(K_n)$;
\item[(3)] $2r(M_{n,j-1})\leqslant q_{n,j}\leqslant r(M_{n,j})$ for $1\leqslant j\leqslant
3b$.
\end{enumerate}
\end{lemma}

The following is an immediately corollary.
\begin{cor}
For any $n\geqslant 1$,
\begin{enumerate}
\item[(1)] $s_{n-1}\leqslant r(L_{n});$

\item[(2)] $r(I_n)\geqslant 2^{3b}r(I_{n-1})$.
\end{enumerate}
\end{cor}

\section{Bounded shape of puzzle pieces}

\setcounter{equation}{0} \noindent We suppose $T(c_0)$ is
persistently recurrent and puzzle pieces $c_0\in
\widetilde{K}_n\subset K_n\subset K_n'$ are constructed as in the
previous section.

For the polynomial case, the following is the key lemma in
\cite{QY}.

\begin{lemma}[\cite{QY}]
$\liminf_{n\to\infty} \mathrm{mod}(K_n'\setminus
\overline{K}_n)>0$.
\end{lemma}

\begin{lemma}There exist a constant $m>0$ depending only on
$b$ and $d_{max}$, and an integer $n_0$ such that
$\mathrm{mod}(K_n'\setminus \overline{K}_n)\geqslant m$ and
$\mathrm{mod}(K_n\setminus \overline{\widetilde{K}_n})\geqslant m$
for all $n\geqslant n_0$.
\end{lemma}
\begin{proof}
In the attracting case, the annuli $K_n'\setminus \overline{K}_n$
and $K_n\setminus \overline{\widetilde{K}_n}$ are always
non-degenerate.

In the parabolic case, there exists an integer $n_0$ such that
$K_n'\setminus \overline{K}_n$ and $K_n\setminus
\overline{\widetilde{K}_n}$ are non-degenerate for $n\geqslant
n_0$. In fact, there is an integer $k_0$ such that
$P_{0}(c_{0})\setminus \overline{P_{k_{0}}(c_{0})}$ is
non-degenerate because $c_0\not\in \cup_{n\geqslant 0}f^{-n}(0)$
and $\cap_{n\geqslant 0}P_n(c_0)=\{c_0\}$. Take $P_{k_{0}}(c_{0})$
as $I_{0}$ in the construction of KSS nest. By Corollary 1,
\begin{align*}
\mathrm{depth}(K_n)-\mathrm{depth}(K_n')&=\mathrm{depth}(\mathcal{A}(I_n))-\mathrm{depth}(\mathcal{B}(I_n))\\
&\geqslant r(I_{n})\rightarrow \infty
\end{align*}
and
\begin{align*}
\mathrm{depth}(\widetilde{K}_n)-\mathrm{depth}(K_n)&=\mathrm{depth}(\mathcal{A}(L_n))-\mathrm{depth}(\mathcal{B}(L_n))\\
&\geqslant r(L_n)\geqslant r(I_n)\rightarrow \infty.
\end{align*}
So there exists an integer $n_{0}$ such that
$$\mathrm{depth}(K_n)-\mathrm{depth}(K_n')\geqslant k_{0}$$
and
$$\mathrm{depth}(\widetilde{K}_n)-\mathrm{depth}(K_n)\geqslant k_{0}$$
for $n\geqslant n_0$. This implies that $K_n'\setminus
\overline{K}_n$ and $K_n\setminus \overline{\widetilde{K}_n}$ are
non-degenerate for $n\geqslant n_0$ because $K_n$ and
$\widetilde{K}_n$ are pullbacks of $I_0=P_{k_{0}}(c_{0})$.

By the same proof of Lemma 5, $\mathrm{mod}(K_n'\setminus
\overline{K}_n)\geqslant \mu$ for some constant $\mu >0$ depending
only on $b$ and $d_{max}$ when $n\geqslant n_0$. See \cite{QY}.

Suppose $f^{b_{n}}(\widetilde{K}_n)= L_{n}$. Let
$h_n=b_{n}+s_{n}+q_{n-1}$ and
$\widetilde{K}_n'=\mathrm{Comp}_{c_{0}}f^{-h_n}(K'_{n-1})$. Since
$$\mathrm{depth}(\widetilde{K}_n)-\mathrm{depth}(K_n)=r_{f^{t_n}(c_0)}(L_n)\geqslant r(L_n)$$
and
\begin{align*}
\mathrm{depth}(\widetilde{K}_n)-\mathrm{depth}(\widetilde{K}_n')&=\mathrm{depth}(K_{n-1})-\mathrm{depth}(K'_{n-1})\\
&=\mathrm{depth}(\mathcal{A}(I_{n-1}))-\mathrm{depth}(\mathcal{B}(I_{n-1}))\\
&\leqslant s_{n-1},
\end{align*}
we conclude that $\widetilde{K}_n'\subset K_n$ from Corollary 1.

From properties (P1), (P2) and (P4),
$$\deg(f^{h_{n}}|_{\widetilde{K}_n'})=\deg(f^{h_n}|_{\widetilde{K}_n})\leqslant D$$
for some $D<\infty$ depending only on $b$ and $d_{max}$. Hence
\begin{align*}
\mod (K_n\setminus \overline{\widetilde{K}}_{n})&\geqslant \mod
(\widetilde{K}_n'\setminus \overline{\widetilde{K}}_{n})\\
&=\frac{\mod (K'_{n-1}\setminus
K_{n-1})}{\deg(f^{h_n}|_{\widetilde{K}_n})}\\
&\geqslant \frac{\mu}{D}.
\end{align*}

Take $m=\frac{\mu}{D}$. This $m$ depends only on $b$ and
$d_{max}$, and satisfied the conditions set out in this lemma.
\end{proof}

\begin{prop} There exists a constant $M_1>0$ such that
$$\mathrm{Shape}(K_{n},c_{0})\leqslant M_1$$
for $n\geqslant n_0$, where $n_0$ is the integer in Lemma 6.
\end{prop}
\begin{proof}
Since $(K'_n\setminus K_n)\cap \mathrm{orb}([c_{0}])
=(K_{n}\setminus \widetilde{K}_{n})\cap \mathrm{orb}([c_{0}])
=\emptyset$ and $f^{p_{n+1}}(K_{n+1})=K_n$, it follows
$f^{p_{n+1}}(c_0)\in \widetilde{K}_{n}$. Let
$\Omega_n'=\mathrm{Comp}_{c_{0}}f^{-p_{n+1}}(K_n')$ and
$\widetilde{\Omega}_n=\mathrm{Comp}_{c_{0}}f^{-p_{n+1}}(\widetilde{K}_{n})$.
Then
$$2\leqslant \deg(f^{p_{n+1}}|_{\Omega_n'})=\deg(f^{p_{n+1}}|_{K_{n+1}})
=\deg(f^{p_{n+1}}|_{\widetilde{\Omega}_n})\leqslant D_1$$ for some
constant $D_1<\infty$ depending only on $b$ and $d_{max}$.

Let $\varphi:(\Delta,V,\widetilde{V})\rightarrow
(K_n',K_n,\widetilde{K}_n)$ and
$\psi:(\Delta,U,\widetilde{U})\rightarrow
(\Omega_n',K_{n+1},\widetilde{\Omega})$ be conformal maps with
$\varphi (0)=c_0$ and $\psi (0)=c_0$.  Let $g=\varphi^{-1}\circ
f^{p_{n+1}}\circ \psi$. Then
$g:(\Delta,U,\widetilde{U})\rightarrow (\Delta,V,\widetilde{V})$
is a properly holomorphic map with
$$\deg(g\mid_{\widetilde{U}})=\deg(g\mid_{U})=\deg(g\mid_{\Delta})=D_n.$$

By Koebe Distortion Theorem and Lemma 2, there exists a constant
$K$ depending only on $b$ and $d_{max}$ such that
$$\mathrm{Shape}(K_{n+1},c_{0})\leqslant K\cdot \mathrm{Shape}(K_{n},c_{0})^{\frac{1}{2}}$$
hold for all $n\geqslant n_0$. We conclude that
$$\mathrm{Shape}(K_n,c_{0})\leqslant M_1$$
for some constant $M_1>0$ when $n\geqslant n_0$.
\end{proof}

\section{Proofs of Theorem 1 and Theorem 2}

\noindent Recall the definition of $\mathrm{Crit}$ in Section 3.
Let $X_1=\cup_{n\geqslant 0}f^{-n}(\textrm{Crit})$ in the
attracting case and
$$X_1=(\cup_{n\geqslant
0}f^{-n}(\textrm{Crit}))\bigcup (\cup_{n\geqslant 0}f^{-n}(0))$$
in the parabolic case, where $0$ is the parabolic fixed point.

For any $x\in J(f)\setminus X_{1}$, let
\begin{eqnarray*}
&&\mathrm{Crit}(x)=\{c\in\mathrm{Crit}\,|\,\,x\to c\}\\
&&\textrm{Crit}_\textrm{n}(x)=\{c\in\mathrm{Crit}(x)\,| \,\, T(c) \textrm { is non-critical}\},\\
&&\textrm{Crit}_\textrm{p}(x)=\{c\in\mathrm{Crit}(x)\,| \,\, T(c) \textrm { is persistently recurrent}\},\\
&&\textrm{Crit}_\textrm{r}(x)=\{c\in\mathrm{Crit}(x)\,|\, \, T(c) \textrm { is reluctantly recurrent}\},\\
&&\textrm{Crit}_\textrm{en}(x)=\{c^\prime\in\mathrm{Crit}(x)\,|\,
\, c^\prime \not\to c^\prime
\textrm{ and } c^\prime \to c \textrm{ for some } c\in\textrm{Crit}_\textrm{n}(x) \},\\
&&\textrm{Crit}_\textrm{ep}(x)=\{c^\prime\in\mathrm{Crit}(x)\,|\,
\, c^\prime \not\to c^\prime
\textrm{ and } c^\prime \to c \textrm{ for some } c\in \textrm{Crit}_\textrm{p}(x)\},\\
&&\textrm{Crit}_\textrm{er}(x)=\{c^\prime\in\mathrm{Crit}(x)\,|\,
\, c^\prime \not\to c^\prime \textrm{ and } c^\prime \to c
\textrm{ for some } c\in \textrm{Crit}_\textrm{r}(x)\}.
\end{eqnarray*}
Then
$$\mathrm{Crit}(x)=\textrm{Crit}_\textrm{n}(x)\cup \textrm{Crit}_\textrm{p}(x)\cup \textrm{Crit}_\textrm{r}(x)\cup
\textrm{Crit}_\textrm{en}(x) \cup \textrm{Crit}_\textrm{ep}(x)\cup
\textrm{Crit}_\textrm{er}(x).$$

Further let
\begin{eqnarray*}
&&X_2=\{x\in J(f)\setminus X_1\,|\,\textrm{Crit}(x)=\emptyset
\textrm{ or }
\textrm{Crit}_\textrm{n}(x)\cup\textrm{Crit}_\textrm{r}(x)\neq\emptyset\},\\
&&X_3=\{x\in J(f)\setminus
X_1\,|\,\textrm{Crit}(x)=\textrm{Crit}_\textrm{p}(x)\cup
\textrm{Crit}_\textrm{ep}(x),\, \textrm{Crit}_\textrm{ep}(x)\neq\emptyset\},\\
&&X_4=\{x\in J(f)\setminus
X_1\,|\,\textrm{Crit}(x)=\textrm{Crit}_\textrm{p}(x)\neq\emptyset\}.
\end{eqnarray*}
Then
$$J(f)=\bigcup_{i=1}^{4}X_{i}.$$

\begin{lemma} For any $x\in X_2\cup X_3$, there exist a puzzle piece $P_0$ of depth $0$ and
infinitely many $i_n$ such that
$$\deg(f^{i_n}:P_{i_n}(x)\rightarrow P_0)\leqslant D$$
for some constant $D<\infty$ depending on $x$.
\end{lemma}
\begin{proof} There are four possibilities.

(1) $T(x)$ is non-critical, i.e. $\textrm{Crit}(x)=\emptyset$.
There exists an integer $n_0\geqslant 0$ such that $(n_0,j)$ is
not a critical position for all $j>0$. For any $n\geqslant 1$,
$\deg(f^n|_{P_{n_0+n}(x)})\leqslant \deg(f|_{P_{n_0+1}(x)})$. The
degrees of these maps $$f^{n_0+n}:P_{n_0+n}(x)\rightarrow
P_0(f^{n_0+n}(x))$$ has an upper bound $D<\infty$. Take a
subsequence $i_n$ of $n_0+n$ such that $P_0(f^{i_n}(x))=P_0$ for
some fixed puzzle piece $P_0$. Then
$$\deg(f^{i_n}:P_{i_n}(x)\rightarrow P_0)\leqslant D$$
for all $n$.

(2) $x\to c$ for some $c\in \textrm{Crit}_\textrm{n}(x)$. From
(1), there are a puzzle piece $P_0$, a positive integer $N_1$ and
infinitely many $j_n$ such that
$$\deg(f^{j_n}:P_{j_n}(c)\rightarrow P_0)\leqslant N_1$$
for all $n$. For each $n$, let $l_n$ be the first moment such that
$f^{l_n}(x)\in P_{j_n}(c)$, $i.e.,\,\,(j_n, l_n)$ is the first
$c$-position on the $j_n$-th row in $T(x)$. By tableau rules (T1)
and (T2), there is at most one $c^\prime$-position on the diagonal
$$\{(k,m)\,|\, \, k+m = j_n+l_n, \quad j_n< k \leqslant j_n+l_n\}$$
for any $c^\prime \in \textrm{Crit}(x) - \{c\}.$ There exists a
positive integer $N_2$ depending on $\textrm{Crit}(x)$ such that
$$\deg(f^{l_n}:P_{j_n+l_n}(x)\rightarrow P_{j_n}(c))\leqslant
N_2.$$ Take $i_n=j_n+l_n$ and $D=N_1+N_2$. Then
$$\deg(f^{i_n}:P_{i_n}(x)\rightarrow P_0)\leqslant D$$
for all $n$.

(3) $x\to c$ for some $c\in \textrm{Crit}_\textrm{r}(x)$. There
exist an integer $n_0\geq 0$, $c^\prime \in [c]$,  $c_1 \in [c]$
and infinitely many integers $k_n\geqslant 1$ such that
$\{P_{n_0+k_n}(c^\prime)\}_{n\geqslant 1}$ are children of
$P_{n_0}(c_1)$. Since $c^\prime \in [c]$, we have $x\to c^\prime$.
For each $n$, let $m_n$ be the first moment such that
$f^{m_n}(x)\in P_{n_0+k_n}(c^\prime)$. There is at most one
$\widetilde{c}$-position on the diagonal
$$\{(n,m)\,|\,\, n+m = n_0+k_n+m_n, \quad n_0+k_n < n \leqslant n_0+k_n+m_n\}$$
in $T(x)$ for any $\widetilde{c} \in \textrm{Crit}(x) -
\{c^\prime\}.$ Therefore, $f^{m_n+k_n}(P_{n_0+k_n+m_n}(x))=
P_{n_0}(c_1)$ and there is an integer $N_3$ independent of $n$
such that
$$\deg (f^{m_n+k_n}|_{P_{n_0+k_n+m_n}(x)})\leqslant N_3
< \infty$$
for any $n\geqslant 1$. Take $i_n=n_0+k_n+m_n$,
$D=N_3+\deg(f^{n_0}|_{P_{n_0}(c_1)})$ and
$P_0=f^{n_0}(P_{n_0}(c_1))$. Then
$$\deg(f^{i_n}:P_{i_n}(x)\rightarrow P_0)\leqslant D$$
for all $n$.

(4) $x\in X_3$, i.e.
$\textrm{Crit}_\textrm{n}(x)\cup\textrm{Crit}_\textrm{r}(x)=\emptyset$
and $\textrm{Crit}_\textrm{ep}(x)\neq\emptyset$. Take
$c_0\in\textrm{Crit}_\textrm{ep}(x)$. Let$\{(0,j_n\}_{n\geqslant 1
}$ be all $c_0$-positions in $T(x)$. We claim that there is at
most one $c$-position on the diagonal
$$\{(n,m)\,|\,\, n+m = j_n, \quad 0< n \leqslant j_n\}$$
for all $c\in \textrm{Crit}(x)$. If this is false, there are at
least two $c$-positions on this diagonal for some $c\in
\textrm{Crit}(x)$. This means $c\in \textrm{Crit}_\textrm{p}(x)$
and $c_0\in F(c)$. By Lemma 1, $F(c)=[c]$ and $c_0\in
\textrm{Crit}_\textrm{p}(x)$. It contradicts with
$c_0\in\textrm{Crit}_\textrm{ep}(x)$. So the above claim is true.
There exists a positive integer $D$ such that
$$\deg(f^{j_n}:P_{j_n}(x)\rightarrow P_0(f^{j_n}(x)))\leqslant D.$$
Take a subsequence $\{i_n\}$ of $\{j_n\}$ such that
$P_0(f^{i_n}(x))=P_0$ for a fixed $P_0$. Then
$$\deg(f^{i_n}:P_{i_n}(x)\rightarrow P_0)\leqslant D$$
for all $n$.
\end{proof}

\begin{prop}
$\mathrm{mes}(X_1\cup X_2\cup X_3)=0$, where $\mathrm{mes}$
denotes the Lebesgue measure on the complex plane $\mathbb{C}$.
\end{prop}
\begin{proof}
It is sufficient to prove that any point $x\in X_2\cup X_3$ is not
a density point of $J(f)$.

From Lemma 7, for any point $x\in X_2\cup X_3$, there exist a
puzzle piece $P_0$ and infinitely many $i_n$ such that
$$\deg(f^{i_n}:P_{i_n}(x)\rightarrow P_0)\leqslant D$$
for some constant $D<\infty$ depending on $x$.

{\bf {The attracting case}}. There is a subsequence of $\{i_n\}$,
say itself, such that $\{f^{i_n}(x)\}$ converges to some point
$x_0\in J(f)$. We assume that $f^{i_n}(x)\in P_1(x_0)$ for all
$n$. It is obvious that there is a disk $D(y,r_0)$ in
$P_1(x_0)\cap F(f)$ for some constant $r_0>0$. By distortion
results for holomorphic $p$-valent mappings(see \cite{CJY},
\cite{Hai},\cite{Sh},\cite{ST} and \cite{Yin}), there are
constants $1\leqslant M <\infty$ and $0<\lambda < 1$ depending on
$x$ such that
$$\text{Shape}(P_{i_n+1}(x),x)\leqslant M$$
and
$$\frac{\textrm{mes}(P_{i_n+1}(x)\cap J(f))}{\textrm{mes}(P_{i_n+1}(x))}\leqslant \lambda$$
for all $n$. Since $\cap_{n\geqslant 0}P_{i_n+1}(x)=\{x\}$, the
point $x$ is not a density point of $J(f)$.

{\bf {The parabolic case}}. If there are a puzzle piece $P_{n_0}$
of depth $n_0$ compactly contained in $P_0$ and infinitely many
$i_n$ such that $f^{i_n}(x)\in P_{n_0}$, we can prove that $x$ is
not a density point of $J(f)$ by the same argument as in the
attracting case. Otherwise, there is a subsequence of
$\{f^{i_n}(x)\}$, say itself, converges to some point $x_0\in
\partial P_0$. The point $x_0$ belongs to $\cup_{n\geqslant
0}f^{-n}(0)$. We assume that $0$ is the parabolic fixed point as
before. In the construction of the Branner-Hubbard puzzle, the
flower petal $U_0$ can be chosen such that there exists a sector
$S\subset P_0$ with the vertex at $x_0$. For all $n$, let $r_n> 0$
be the distance from the point $f^{i_n}(x)$ to the boundary of the
sector $S$. Let
$$\widetilde{\Omega}_n(x)=\textrm{Comp}_x(f^{-i_n}(D(f^{i_n}(x),r_n)))$$
and
$$\Omega_n(x)=\textrm{Comp}_x(f^{-i_n}(D(f^{i_n}(x),\frac{1}{2}r_n))).$$
Then
$$\deg(f^{i_n}:\widetilde{\Omega}_n(x)\rightarrow D(f^{i_n}(x),r_n))\leqslant D.$$
By the Leau-Fatou Flower Theorem, there is a disk
$$D(y_n,\frac{1}{4}r_n)\subset D(f^{i_n}(x),\frac{1}{2}r_n)\cap F(f)$$
for large $n$. By the same argument as above, $x$ is not a density
point of $J(f)$.
\end{proof}

Recall the definition of invariant line fields in Section 1.

Let $\mathcal{H}(f)$ be the collection of all holomorphic maps
$h:U\rightarrow V$, where $U$, $V$ are open sets such that there
exist $i,\textrm{ }j\in \mathbb{N}$ with $f^i\circ h=f^j$ on $U$.

The following proposition due to Weixiao Shen is a criterion to
test whether a rational map carries invariant line field on its
Julia set or not.

\begin{prop}[\cite{Sh}]
Let $f$ be a rational map of degree $\geqslant 2$ and $x$ be a
point in $J(f)$. If there are a constant $C>1$, a positive integer
$N\geqslant 2$ and a sequence $h_n:U_n\rightarrow V_n$ in
$\mathcal{H}(f)$ with the following properties:
\begin{enumerate}
\item[(1)] $U_n$, $V_n$ are topological disks and
$$\mathrm{diam}(U_n)\rightarrow 0\textrm{ and } \mathrm{diam}(V_n)\rightarrow 0$$
as $n\rightarrow\infty$.
\item[(2)] $h_{n}$ is a proper map of degree between $2$ and
$N$.
\item[(3)] For some $u\in U_n$ such that $h_n'(u)=0$ and
for $v=h_n(u)$ we have $\textrm{Shape}(U_n,u)\leqslant C$ and
$\textrm{Shape}(V_n,v)\leqslant C$.
\item[(4)] $d(U_n,x)\leqslant C\cdot
\mathrm{diam}(U_n)\textrm{ and } d(V_n,x)\leqslant C\cdot
\mathrm{diam}(V_n).$
\end{enumerate}
Then for any $f$-invariant line field $\mu$, $\mu(x)=0$ or $\mu$
is not almost continuous at $x$.
\end{prop}

\begin{prop}
Suppose $\mu$ is an invariant line field on the Julia set $J(f)$.
If $x\in X_4$, then $\mu (x)=0$ or $\mu$ is not almost continuous
at $x$.
\end{prop}
\begin{proof}
For $x\in X_4$,
$\textrm{Crit}(x)=\textrm{Crit}_\textrm{p}(x)\neq\emptyset$. Let
$c_0\in \textrm{Crit}_\textrm{p}(x)$ and let
$\widetilde{K}_n\subset K_n\subset K_n'$ be puzzle pieces around
$c_0$ as in Section 3.

For any $n\geqslant n_0$, the annuli $K_n'\setminus K_n$ and
$K_{n}\setminus \widetilde{K}_{n}$ are non-degenerate. Let $l_n$
be the first moment such that $f^{l_n}(x)\in \widetilde{K}_n$. Let
$\widetilde{V}_n(x)=\textrm{Comp}_x(f^{-l_n}(\widetilde{K}_n))$,
$V_n(x)=\textrm{Comp}_x(f^{-l_n}(K_n))$ and
$V_n'(x)=\textrm{Comp}_x(f^{-l_n}(K_n'))$. For large $n$, the
puzzle piece $V_n'(x)$ contains no critical points. Let $v_n>0$ be
the smallest integer such that $f^{v_n}(\widetilde{V}_n(x))$
contains a critical point $c\in [c_0]$. Set
$\widetilde{\Lambda}_n=f^{v_n}(\widetilde{V}_n(x))$,
$\Lambda_n=f^{v_n}(V_n(x))$ and $\Lambda_n'=f^{v_n}(V_n'(x))$. See
Figure~\ref{fig3}.

\begin{figure}[h]
\psfrag{Tx}{$T(x):$}

\psfrag{c0}{$c_{0}$} \psfrag{c0t}{$c$}

\psfrag{Vnxp}{$V_n'(x)$} \psfrag{Vnx}{$V_n(x)$}

\psfrag{Vnxt}{$\widetilde{V}_n(x)$} \psfrag{lambdap}{$\Lambda'$}

\psfrag{lambda}{$\Lambda$} \psfrag{lambdat}{$\widetilde{\Lambda}$}

\psfrag{Knp}{$K_n'$} \psfrag{Kn}{$K_n$}

\psfrag{Knt}{$\widetilde{K}_n$} \psfrag{fvn}{$f^{v_n}$} \centering
\includegraphics[width=10cm]{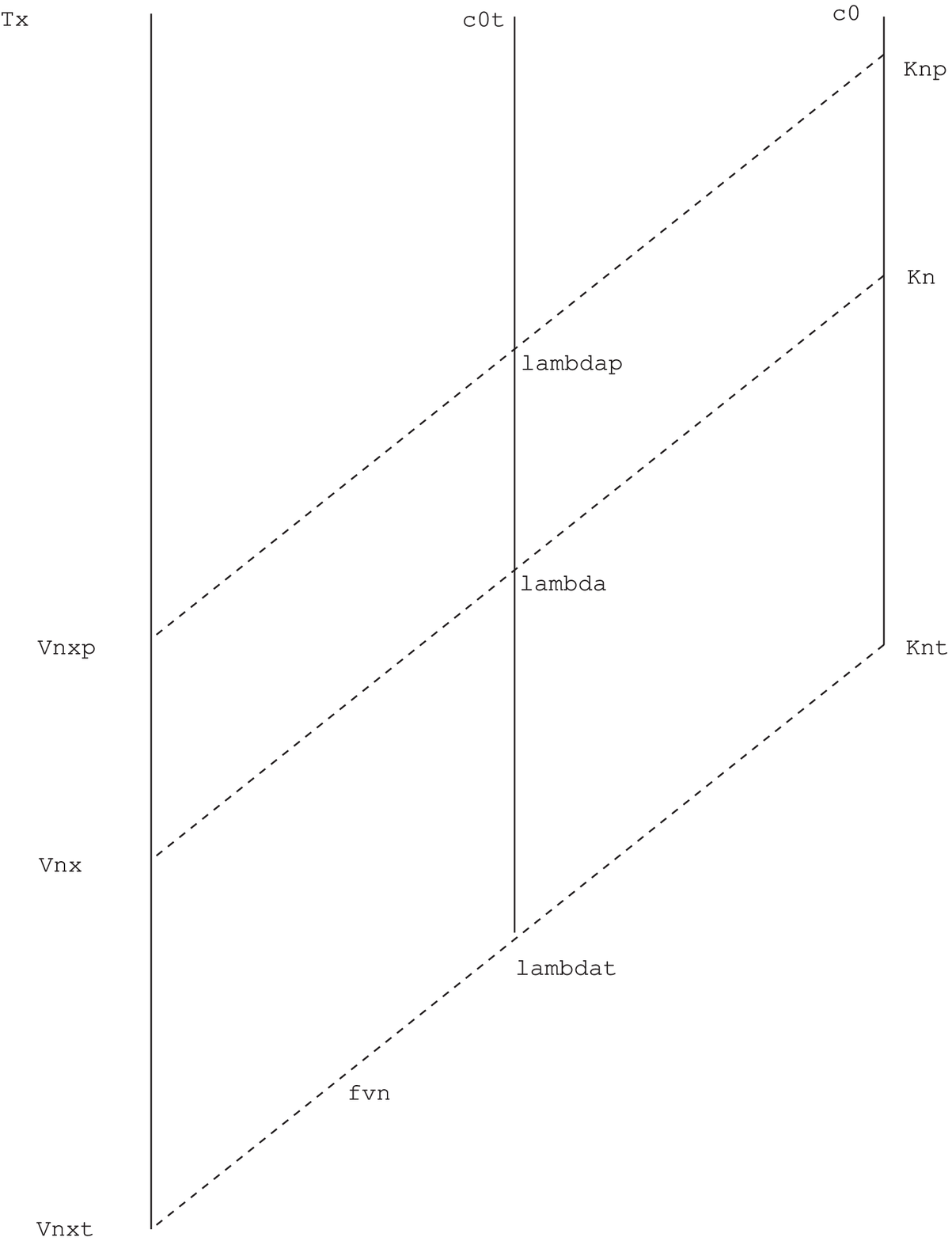} \caption{ } \label{fig3}
\end{figure}

From the conditions
$$(K'_n\setminus K_n)\cap \mathrm{orb}([c_{0}])=(K_{n}\setminus
\widetilde{K}_{n})\cap \mathrm{orb}([c_{0}])=\emptyset$$ and
$\textrm{Crit}(x)=\textrm{Crit}_\textrm{p}(x)$, we know that
$f^{v_n}:V_n'(x)\rightarrow \Lambda_n'$ is conformal and
$$2\leqslant \deg(f^{l_n}|_{V_n'(x)})=\deg(f^{l_n}|_{V_n(x)})=\deg(f^{l_n}|_{\widetilde{V}_n(x)})\leqslant D_2$$
for some constant $D_2$ depending only on $\textrm{Crit}(x)$. By
Proposition 1 and Lemma 2,
$$\textrm{Shape}(V_n(x),x)\leqslant M_2$$
and
$$\textrm{Shape}(\Lambda_n,c)\leqslant M_2$$
for some constant $M_2>0$.

Let $\widetilde{\Gamma}_n=\mathcal{L}_c(\widetilde{\Lambda}_n)$,
$\widetilde{U}_n=\mathcal{L}_x(\widetilde{\Gamma}_n)$ and
$f^{u_n}(\widetilde{U}_n)=\widetilde{\Lambda}_n$. Further let
$\widetilde{U}_n(x)=\textrm{Comp}_xf^{-u_n}(\widetilde{\Lambda}_n)$,
$U_n(x)=\textrm{Comp}_xf^{-u_n}(\Lambda_n)$ and
$U_n'(x)=\textrm{Comp}_xf^{-u_{n}}(\Lambda_n')$. Then
$$2\leqslant \deg(f^{u_n}|_{U_n'(x)})=\deg(f^{u_n}|_{U_n(x)})=\deg(f^{u_n}|_{\widetilde{U}_n(x)})\leqslant D_3$$
for some constant $D_3$ depending only on $\textrm{Crit}$. There
exists a positive constant $M_3>0$ such that
$$\textrm{Shape}(U_n,x)\leqslant M_3.$$

For each large $n$, define $h_n=f^{-v_n}\circ f^{u_n}$. Then
$h_n:U_n(x)\rightarrow V_n(x)$ is a properly holomorphic mapping
of degree between $2$ and some constant $N$. All conditions in
Proposition 3 are satisfied, hence the lemma holds.
\end{proof}

\begin{proof}[Proof of Theorem 1] If $f$ has an invariant line
field $\mu$ on the Julia set $J(f)$, then there exists a positive
measure subset $E$ of $J(f)$ such that $\textrm{support}(\mu)=E$.
Since $\mu$ is measurable in $\mathbb{C}$, almost every point in
$\mathbb{C}$ is almost continuous. From Proposition 4 and
Proposition 2, $\textrm{mes}(X_4\cap \textrm{support}(\mu))=0$
$\textrm{mes}(\textrm{support}(\mu))=0$. It is a contradiction. So
$f$ carries no invariant line field on its Julia set.
\end{proof}

Recall that the Teichm\"{u}ller space of a rational map $f$ of
degree $d$ is defined by
$$\mathrm{Teich}(\hat{\mathbb{C}},f)=\{g\in \mathrm{Rat}_{d}\,|\,g \textrm{ is quasiconformally conjugated with }
f\}/\mathrm{Aut}(\hat{\mathbb{C}}).$$

C. McMullen and D. Sullivan have given a formula to compute the
dimension of $\mathrm{Teich}(\hat{\mathbb{C}},f)$ in \cite{McS}.

\begin{theoremC}[\cite{McS}] The dimension of the Teichm\"{u}ller
space of a rational map is given by
$N=N_{AC}+N_{HR}+N_{LF}-N_{P}$, in which

$\bullet$ $N_{AC}$ is the number of foliated equivalence classes
of acyclic critical point in the Fatou set;

$\bullet$ $N_{HR}$ is the number of cycles of Herman rings;

$\bullet$ $N_{LF}$ is the number of ergodic line field on the
Julia set;

$\bullet$ $N_{P}$ is the number of the parabolic cycles.
\end{theoremC}

\begin{proof}[Proof of Theorem 2]
Let $f$ be a structurally stable rational map with a Cantor Julia
set. From Theorem B in Section 1,
$\dim(\mathrm{Teich}(\hat{\mathbb{C}},f))=2d-2$.

By the implicit function theorem, any rational map with
indifferent cycles is structurally unstable, therefore $f$ has no
Siegel disks or parabolic basins. Moreover, any rational map with
Herman rings is also structurally unstable by a result due to
Ma\~{n}\'{e}, see \cite{Ma}. These imply that $N_{HR}=N_{P}=0$.
From Theorem 1, $N_{LF}=0$. We conclude that $N_{AC}=2d-2$ and all
critical points are in the attracting Fatou components. It means
that $f$ is hyperbolic.
\end{proof}

\begin{acknowledgements} This paper is based on a joint
work of the first author with Weiyuan Qiu and works of Weixiao
Shen. We are grateful to them for their helpful ideas and
discussions.

This work was partially supported by the National Natural Science
Foundation of China.
\end{acknowledgements}

\bibliographystyle{amsplain}

\noindent{Yongcheng Yin}

{School of Mathematical Sciences}

{Fudan University}

{Shanghai, 200433}

{P.R.China}

{yin@zju.edu.cn}

\noindent{Yu Zhai}

{Department of Mathematics}

{Zhejiang University}

{Hangzhou, 310027}

{P.R.China}

{dyu@zju.edu.cn}
\end{document}